\documentclass[11pt]{amsart}

\usepackage{hyperref}
\usepackage{xcolor}
\hypersetup{colorlinks=true,urlcolor=black,linkcolor=black,citecolor=black}
\usepackage{amscd,amsmath,amssymb}
\usepackage[mathscr]{euscript}
\usepackage[all]{xy}
\usepackage{fullpage}

\usepackage[charter]{mathdesign}
\usepackage[T1]{fontenc}

\usepackage{cite}

\makeatletter \@addtoreset{equation}{section}\makeatother

\DeclareMathOperator*{\krest}{\otimes}

\newtheorem{theorem}{Theorem}[section]

\newtheorem{proposition}[theorem]{Proposition}
\newtheorem{corollary}[theorem]{Corollary}

\newcommand{\ZZ}{{\mathbb{Z}/2}}

\newcommand{\obj}{{\mathfrak{a}}}

\newcommand{\Al}{{\mathbb{A}^{\!1}}}
\newcommand{\Am}{{\mathbb{A}^{\!m}}}
\newcommand{\Lmo}{{{\Lambda}_{m}}}
\newcommand{\Lm}{{{\Lambda}^{\rm \!ext}_{m}}}
\newcommand{\Lo}{{{\Lambda}_{0}}}
\newcommand{\A}{\C[\ttt]}
\newcommand{\C}{{\mathbb{K}}}
\newcommand{\CC}{{\mathbb{C}}}
\newcommand{\ev}{\mathsf{ev}}
\newcommand{\od}{\mathsf{od}}

\newcommand{\rrC}{{\mathrm{C}}}
\newcommand{\rrCn}{{\mathrm{C}}}
\newcommand{\rC}{{{\mathrm{C}}}}
\newcommand{\rCu}{{{\mathscr{C}}}}
\newcommand{\rrCu}{{{\mathscr{C}}}}
\newcommand{\rZ}{{\mathsf{Z}}}
\newcommand{\rg}{{\mathrm{R}\Gamma}}
\newcommand{\bb}{{b}}
\newcommand{\B}{{{B}}}
\newcommand{\BB}{{{\beta}}}

\newcommand{\ee}{{{e}}}

\newcommand{\EE}{{{E}}}

\newcommand{\cA}{{\mathfrak{A}}}

\newcommand{\cPhi}{{{\Phi}}}
\newcommand{\cP}{{\mathfrak{perf}}}
\newcommand{\cK}{{\mathscr{K}}}
\newcommand{\cX}{{{X}}}

\newcommand{\cE}{{\mathscr{E}}}
\newcommand{\cO}{{\mathscr{O}}}
\newcommand{\cD}{{\mathcal{D}}}
\newcommand{\cF}{{\mathscr{F}}}
\newcommand{\cI}{{\mathscr{I}}}
\newcommand{\rM}{{\rC}}
\newcommand{\cM}{{\rCu}}
\newcommand{\rCP}{{\mathrm{Tot}_{\yy,u}}}
\newcommand{\rHP}{{\mathrm{HP}}}

\newcommand{\rH}{{\mathrm{H}}}
\newcommand{\uu}{{(\!(\!u\!)\!)^\mathsf{gr}}}
\newcommand{\ttt}{\boldsymbol{\!t\!}}
\newcommand{\ww}{\boldsymbol{f}}
\newcommand{\yy}{\boldsymbol{\tau}}
\newcommand{\yyy}{{[\!\yy\!]}}
\newcommand{\tc}{{\text{\guillemotleft}t\cdot\text{\guillemotright}}}
\newcommand{\id}{{{\rm id}}}

\newcommand{\Dm}{{A_m}}
\newcommand{\DNC}{{\widetilde{A}_m}}

\allowdisplaybreaks

\title{Exponentially twisted cyclic homology}
\author{Dmytro Shklyarov}

\email{dmytro.shklyarov@mathematik.tu-chemnitz.de}

\begin{document}
\begin{abstract} By a theorem of Bernhard Keller the de Rham cohomology of a smooth variety is isomorphic to the periodic cyclic homology of the differential graded category of perfect complexes on the variety. Both the de Rham cohomology and the cyclic homology can be twisted by the exponential of a regular function on the variety. We explain that the isomorphism holds true in the twisted setting and draw some conclusions on derived invariance of the algebraic Gauss-Manin systems associated with regular functions.
\end{abstract}

\maketitle

\baselineskip 1.4pc

\section{Introduction}  

According to  \cite{Gr}, the cohomology $\rH^{*}(\cX^{\rm \!an}, \CC)$ of the analytic variety $\cX^{\rm \!an}$ associated with a smooth complex variety\footnote{In this note ``variety'' means ``quasi-projective variety''.} $\cX$ can be computed algebraically:
\begin{equation}\label{gr}
\rH^k(\cX^{\rm \!an})\simeq \rH^k_{\rm DR}(\cX):=\mathbb{H}^k(\Omega^*_\cX, d).
\end{equation}
Here, and in the rest of the paper, $(\Omega^*_\cX, d)$ denotes the algebraic de Rham complex and $\mathbb{H}^*$ denotes the hypercohomology in the Zariski topology. 
A more recent result \cite{K0} gives the right-hand side of (\ref{gr}) a categorical interpretation:
\begin{equation}\label{ke}
\bigoplus\limits_{k\,\text{\rm even/odd}} \rH^{k}_{\rm DR}(\cX)\simeq\rHP_{\ev/\od}(\cP\,\cX,\cP_\circ\cX)
\end{equation}
where $\cP\,\cX$ is the dg category of perfect complexes on $\cX$,  $\cP_\circ\cX$ is the dg subcategory of acyclic perfect complexes, and $\rHP_*$, $*\in\ZZ$, denotes the periodic cyclic homology.
A remarkable consequence of these theorems is that $I(\cX):=\bigoplus_{k\,\text{\rm even/odd}} \rH^{k}(\cX^{\rm \!an})$ turns out to be a derived invariant of smooth varieties:  $I(\cX)\simeq I(\cX')$ for any smooth Fourier-Mukai partners $\cX$ and $\cX'$. 

This note grew out of an attempt to generalize these results to the setting of pairs $(\cX, f)$ where $\cX$ is a smooth variety and $f$ is a regular function on it. 
The isomorphism (\ref{gr})  has a counterpart in this world, namely \cite{DS,Sab}
\begin{equation}\label{b}
\rH^k(\cX^{\rm \!an}, f^{-1}(\!t\!))\simeq \rH^k_{\rm DR}(\cX, f):=\mathbb{H}^k(\Omega^*_\cX, d-df) \quad ({\rm Re}\, t\to +\infty)
\end{equation}
Since, informally, $d-df=e^fde^{-f}$, the cohomology in the right-hand side of  (\ref{b}) is called the (exponentially) twisted de Rham cohomology. 
A natural question is whether there is an analog of (\ref{ke}).  

We are actually interested in a more sophisticated version of the twisted de Rham cohomology. Namely, we would like to have a categorical interpretation of the following one-parameter family of twisted cohomology groups\footnote{Throughout the text, $V[\!x\!]$, $V[\!x^{\pm1}\!]$, $V(\!(\!x\!)\!)$ stand for the spaces of polynomials resp. Laurent series resp. formal Laurent series in $x$ with coefficients in $V$.}$^,$\footnote{In the main body of the paper we discuss a ``multi-parameter'' generalization of (\ref{gm1}) -- a $\cD_\Am$-module associated with a variety $\cX$ and a collection $\ww=(f_1,\ldots, f_m)$ of functions on it \cite{AS,DMSS} -- but for expository purposes we are focusing on the one-parameter case here.} 
\begin{eqnarray}\label{gm1}
\mathbb{H}^*(\Omega^*_\cX[\!\tau\!],\,  d-\tau df)\quad (\text{informally,}\,\,\,\, d-\tau df=e^{\tau f}de^{-\tau f})
\end{eqnarray} 
{\it together} with its $\cD_\Al$-module structure induced by the endomorphisms 
\[
\tau:=\tau\cdot,\quad \partial_\tau:=\frac{\partial}{\partial\tau}-f \,\,(=e^{\tau f}\frac{\partial}{\partial\tau}\,e^{-\tau f})
\]
of the complex $(\Omega^*_\cX[\tau], d-\tau df)$.  This $\cD_\Al$-module, known as the Laplace (transform of the) Gauss-Manin system of $f$ \cite{Sa0}, and various versions of it, e.~g.
\begin{eqnarray}\label{gm2}
\mathbb{H}^*(\Omega^*_\cX[\!\tau^{\pm1}\!], d-\tau df),\quad \mathbb{H}^*(\Omega^*_\cX(\!(\!\tau^{-1}\!)\!), d-\tau df),
\end{eqnarray}
encode a great deal of information on the topology of the map $f:\cX^{\rm \!an}\to\CC$ and form the basis of the Hodge theory of the pairs $(\cX,f)$; see e.~g.\cite{HSe,KKP,KS,Sa0,Sa,Sa1,SaSa}.

Thus, we are looking for a version of (\ref{ke}) with (\ref{gm1})  on the left-hand side and a $\cD_\Al$-module of categorical origin on the right-hand side.

The multiplication operator
$
\Omega^*_\cX\stackrel{f\cdot}\to \Omega^*_\cX
$
which enters the above definitions has a well-known analog in the framework of non-commutative differential calculus \cite{Ts}. Let $A$ be an associative unital algebra. Recall \cite{L} that the complex computing $\rHP_*(A)$  is defined as 
$
(\rrC_*(A)\uu, b+uB)
$
where $(\rrC_*(A):=A\otimes (A/\CC)^{\otimes (-*)}, b)$ is the Hochschild complex, 
$B$ is the cyclic differential, 
$u$ is a degree 2 variable, and $\rrC_*(A)\uu$ denotes the subspace in $\rrC_*(A)(\!(\!u\!)\!)$ spanned by the homogeneous formal Laurent series. Given a central element $t\in A$ the analog of $f\cdot$ is the endomorphism of $\rrC_*(A)\uu$ given by the formula
\begin{eqnarray*}
\tc: a_0\otimes a_n\otimes.. \otimes a_1&\mapsto&ta_0\otimes a_n\otimes.. \otimes a_1+\\
&+&u\sum_{p=1}^{n+1}\sum_{l=0}^{p-1}(-1)^{p+nl+1}1\otimes a_l\otimes a_{l-1}\otimes .. \otimes a_0\otimes a_p \otimes t\otimes a_{p-1}\otimes ..\otimes a_{l+1}.
\end{eqnarray*}
Thus, we have a natural non-commutative analog of the $\cD_\Al$-module (\ref{gm1}):
\begin{eqnarray*}
\rHP_*(A,t):=\rH^*(\rrC_*(A)[\tau]\uu, e^{\tau\,\tc}(b+uB)e^{-\tau\,\tc}),\quad \partial_\tau:=e^{\tau\,\tc}\frac{\partial}{\partial\tau}e^{-\tau\,\tc}.
\end{eqnarray*}
This construction extends to the setting of pairs $(\cA, t)$, where $\cA$ is a dg category and $t$ is a closed degree 0 natural transformation of the identity endofunctor of $\cA$, and further to the setting of triples $(\cA,\cA_\circ, t)$ where $(\cA, t)$ are as before and $\cA_\circ$ is a full dg category of $\cA$. Notice that  the function $f$ on $\cX$ defines an endomorphism of every complex of $\cO_\cX$-modules and these endomorphisms altogether form a natural transformation of the identity endofunctor of $\cP\,\cX$. Thus, the triple $(\cP\,\cX,\cP_\circ\cX, f)$ gives rise to a $\cD_\Al$-module $\rHP_*(\cP\,\cX,\cP_\circ\cX, f)$. Our generalization of (\ref{ke}) can be formulated as follows (Theorem \ref{mt1}): One has an isomorphism of $\cD_\Al$-modules 
\begin{eqnarray}\label{gm3}
\bigoplus_{k\,\text{\rm even/odd}} \mathbb{H}^{k}(\Omega^*_\cX[\!\tau\!],\, d-\tau df)\simeq\rHP_{\ev/\od}(\cP\,\cX,\cP_\circ\cX, f).
\end{eqnarray}

It follows from (\ref{gm3}) (cf. Corollary \ref{cor}) that two pairs $(\cX, f)$ and $(\cX', f')$ give rise to isomorphic total even/odd Laplace Gauss-Manin systems if $\cX$ and $\cX'$ are {\it relative} Fourier-Mukai partners, i.~e. if there exists a Fourier-Mukai equivalence between the bounded derived  categories of coherent sheaves $D^b(\cX)$ and $D^b(\cX')$ with kernel in $D^b({\cX\!\times_{\!\Al}\!\cX'})$. Here is an example of how this observation may be utilized. Let $\mathit{\Gamma}$ be a finite subgroup of $SL_n(\CC)$, $Y=\CC^n/\mathit{\Gamma}$, and  $\cX\stackrel{\phi}\to Y\stackrel{\phi'}\gets\cX'$ be two crepant resolutions. Let also  $g$ be a regular function on $Y$ and $f:=g\circ\phi$ and $f':=g\circ\phi'$ be the corresponding pullbacks. The derived Mckay correspondence conjecture asserts that $\cX$ and $\cX'$ are Fourier-Mukai partners. Moreover, it is expected \cite{BP} that there exists a Fourier-Mukai equivalence between $D^b(\cX)$ and $D^b(\cX')$ with kernel in $D^b({\cX\!\times_{\!Y}\!\cX'})$. Thus, in the cases where this assertion is known to hold one may conclude that the even/odd Laplace Gauss-Manin systems associated with $(\cX,f)$ and $(\cX',f')$ are isomorphic because an equivalence with kernel in $D^b({\cX\!\times_{\!Y}\!\cX'})$ can also be viewed, using the embedding  ${\cX\!\times_{\!Y}\!\cX'}\to{\cX\!\times_{\!\Al}\!\cX'}$, as an equivalence with kernel in $D^b({\cX\!\times_{\!\Al}\!\cX'})$.   

Arguments of this sort may prove useful in the study of some instances of mirror symmetry featuring pairs $(\cX,f)$ -- called in this context Landau-Ginzburg models -- as one of the mirror partners. For example, the construction of Landau-Ginzburg models in \cite{AAEKO,E,GKR,KKOY,Sei} involves arbitrary choices (including resolutions of singularities), and one then has to check that invariants of $(\cX,f)$ of interest in mirror symmetry  are independent of the choices. One such invariant is the triangulated category of singularities (of the singular fibers of $f$). The proof of independence of this category of the choices also relies on the notion of relative Fourier-Mukai equivalence; see \cite{BP} and \cite[Sect. 7]{KKOY}. The twisted de Rham cohomology together with its $\cD$-module structure is another invariant of interest \cite{GKR,KKP} but, to the best of our knowledge, its derived invariance was not discussed previously. 

Let us conclude the Introduction by making two more comments. 

Firstly, observe that the right-hand side of (\ref{gm3}) makes sense for a much broader class of spaces than that of smooth varieties and, thus, provides a natural  generalization of the Laplace Gauss-Manin system beyond the conventional smooth setting. (Note, however, that it is only a $\ZZ$-graded generalization.)

Secondly, one can easily get an analog of (\ref{gm3}) for the $\cD$-modules (\ref{gm2}), as will hopefully be clear from the proof of (\ref{gm3}). In this regard, we would like to mention that the second $\cD$-module in (\ref{gm2})  admits a different categorical interpretation. Namely, as demonstrated in \cite{Ef} (cf. also \cite{Shk1}) it is isomorphic (up to a minor detail) to the ordinary periodic cyclic homology of the dg category of  matrix factorizations associated with $(\cX,f)$, a dg enhancement of the aforementioned category of singularities. This approach, however, does not allow one to recover the $\cD$-module (\ref{gm1}).

\noindent{\bf Conventions.}  $\C$ denotes a field of characteristic 0, our ground field. By a {\it variety} we understand a quasi-projective variety over $\C$. All our complexes are cohomological. We will denote the differential in a complex by $d_{?}$ where $?$ is the notation for the complex; the only exceptions are $d$ (without any subscript) and $b$ -- the  de Rham and the Hochschild differentials, respectively.

\noindent{\bf Acknowledgements.} 
This research was supported by a fellowship from the Freiburg Institute for Advanced Studies (FRIAS), Freiburg, Germany.

\section{DG categories over algebras of polynomials}
\subsection{}\label{dgac} The main objects of interest in the present work are the dg categories over the algebra $\C[\ttt]$ of polynomials in $m$ variables\footnote{Throughout the paper, $m$ denotes a fixed non-negative integer. $m=0$ means ``no variables'', i.~e. $\C[\ttt]=\C$.} $\ttt=(t_1,\ldots,t_m)$, {i.~e.} the small categories whose Hom-sets are  complexes of $\C[\ttt]$-modules satisfying the standard axioms. Alternatively, one may treat such categories as dg categories over $\C$ (henceforth, the latter will be referred to simply as dg categories) endowed with an extra structure -- the ``$t_i$-actions''. It is this viewpoint that we will take, so let us formulate it more precisely.

Let $\cA$ be a dg category and $\rZ(\cA)$ denote the commutative algebra of the closed degree 0 natural transformations of the identity endofunctor of $\cA$,  { i.~e.} the elements $t=\{t_\obj\}$ of $\prod\limits_{\obj\in\cA} \cA^0(\obj,\obj)$ satisfying
\begin{eqnarray*}
d_{\cA}(t_\obj)=0,\qquad t_{\obj'}\,a=a\,t_{\obj}
\end{eqnarray*}
for all $\obj, \obj'\in\cA$ and $a\in \cA(\obj,\obj')$  ($d_{\cA}$ stands for the differential on the Hom-complexes of $\cA$). We will call the pairs $(\cA,\ttt)$, where $\cA$ is a dg category and $\ttt=(t_1,\ldots,t_m)$ with $t_i\in\rZ(\cA)$, {\it dg $\A$-categories}. Also by a {\it dg $\A$-functor} from $(\cA,\ttt)$ to $(\cA',\ttt')$ we will understand a dg functor $\cPhi:\cA\to\cA'$ respecting the natural transformations: 
$
\cPhi(t_{i,\,\obj})=t'_{i,\,\cPhi(\obj)}
$, for all $i=1,\ldots,m$.

It is easy to see that ``dg category over $\C[\ttt]$'' and ``dg $\A$-category'' are equivalent notions. The reason we are being cautious with regards to terminology is that the cyclic-like homology of dg $\A$-categories we will introduce later on is defined in terms of the underlying dg categories and is different from the ordinary ``relative'' cyclic homology.

We will also need a more general framework, namely, that of triples $(\cA,\cA_\circ,\ttt)$ where $(\cA,\ttt)$ is as before  and $\cA_\circ$ is a full dg subcategory of $\cA$. We will call these {\it dg $\A$-pairs}. The functors will be required to preserve the subcategories. The reason we need the more general framework is that we want to include the so-called  localization pairs \cite{K0,K1} in the picture. When $(\cA,\cA_\circ)$ is a localization pair, we will call $(\cA,\cA_\circ,\ttt)$ a {\it localization $\A$-pair}.

\subsection{}\label{ex} A natural source of   dg $\A$-categories are varieties over $\Am$.  
Let $\cX$ be a variety and $\ww=(f_1,\ldots,f_m)$ be a collection of $m$ regular functions on $\cX$. We will denote by $\cP(\cX,\ww)$ the localization $\A$-pair $(\cP\,\cX, \cP_{\circ}\cX, \ww)$ where  $\cP\,\cX$ is the dg category of perfect complexes on $\cX$, $\cP_{\circ}\cX$ is its full subcategory of acyclic complexes, and $f_i$ is viewed as defining an element of $\rZ(\cP\,\cX)$ via the component-wise multiplication $\cE^*\stackrel{f_i\cdot }\to \cE^*$, $\cE=\cE^*\in \cP\,\cX$. 

There are a few variations of this construction which should be thought of as being equivalent to the original one but may be more convenient for technical reasons (for instance, for the purpose of constructing explicit dg $\A$-functors). Namely, $(\cP\,\cX,\cP_{\circ}\cX)$ can be replaced by
\begin{itemize}
\item $(\cP^{\rm \,str}\cX,\cP^{\rm \,str}_\circ\cX)$, the strictly perfect complexes (i.~e. the bounded complexes of algebraic vector bundles) and the acyclic strictly perfect complexes;
\item $(\cP^{\rm \,flab}\cX, \cP^{\rm \,flab}_\circ\cX)$, the perfect bounded below complexes of flabby $\cO_\cX$-modules and its subcategory of all the acyclic complexes;
\item $(\cP^{\rm \,inj}\cX, \cP^{\rm \,inj}_\circ\cX)$, the perfect bounded below complexes of injective $\cO_\cX$-modules and its subcategory of all the acyclic complexes;
\item $\cP^{\rm \,inj}\cX$ by itself (which can also be viewed as the pair $(\cP^{\rm \,inj}, 0)$). 
\end{itemize}
Also, one may take any classical generator $\cE$ of the triangulated category $\rH^0(\cP^{\rm \,inj}\cX)$ (here $\rH^0(-)$ denotes the homotopy category) and consider the dg $\A$-algebra (= a dg $\A$-category with a single object) $({\rm End}(\cE), \ww)$ where ${\rm End}(\cE)$ is the dg algebra of endomorphisms of $\cE$ in $\cP^{\rm \,inj}\cX$ and $\ww$ is viewed as a collection of elements of $\rZ({\rm End}(\cE)\!)$.

A large supply of dg $\A$-functors between the above-described dg $\A$-categories is provided by the following construction. Consider two pairs $(\cX,\ww)$ and $(\cX',\ww')$ as above. Let us assume that both $\cX$ and $\cX'$ are smooth. Then, as explained for example in \cite{BP,KKOY}, every complex $\cK\in D^b({\cX\!\times_{\!\Am}\!\cX'})$ with support proper over $\cX'$ gives rise to a relative Fourier-Mukai functor
\begin{equation}\label{phi}
{\it \Phi}_\cK(-)={\rm R}\pi'_*({\rm L}\pi^*(-)\krest\limits_{\cO_{\cX\!\times_{\!\Am}\!\cX'}}^{\rm L}\cK): D^b(\cX)\to D^b(\cX')
\end{equation}
where $\pi$ and $\pi'$ are the natural morphisms ${\cX\!\times_{\!\Am}\!\cX'}\to\cX$ and ${\cX\!\times_{\!\Am}\!\cX'}\to\cX'$, respectively. Let us pick a bounded below resolution $\cK\to \cI$ with injective components. Then
\begin{equation}\label{Phi}
\cPhi_{\cI}(-)=\pi'_*(\pi^*(-)\krest\limits_{\cO_{\cX\!\times_{\!\Am}\!\cX'}}\cI)
\end{equation}
can be shown to induce a dg functor of localization pairs $(\cP^{\rm \,str}\cX, \cP^{\rm \,str}_{\circ}\cX)\to (\cP^{\rm \,flab}\cX', \cP^{\rm \,flab}_{\circ}\cX')$. Since, by definition, $\pi^*(f_i)=\pi'^*(f'_i)$, this dg functor is a dg $\A$-functor $\cP^{\rm \,str}(\cX,\ww)\to\cP^{\rm \,flab}(\cX',\ww')$.

\section{Twisted cyclic homology of dg $\A$-categories} 
\subsection{}\label{mc} Let us start by introducing some terminology and notation.

Let  $\Lmo$ denote the (super)commutative unital algebra generated by $\BB_{i}$  ($i=0,\ldots, m$), each of degree $-1$.
We will view $\Lmo$ as a dg algebra with the trivial differential.
The (left) dg $\Lmo$-modules can be viewed as a generalization of the mixed complexes  \cite{Ka}; the latter correspond to the case $m=0$. 

Let us fix -- once and for all -- a degree 2 variable $u$ and  a collection of degree 0 variables $\yy=(\tau_1,\ldots,\tau_m)$ and consider the following functor from the category of dg $\Lmo$-modules to the category of complexes of $\C\yyy$-modules: 
\[
\rC=(\rrC^*, d_\rC)\mapsto\rCP(\rC):=(\rrC^*\yyy\uu, \, d_\rC+u\BB_0+u(\tau_1\BB_1+\ldots+\tau_m\BB_m)\!).
\] 
Here, we recall, $\rrC^*\yyy\uu$ denotes the subspace in $\rrC^*\yyy(\!(\!u\!)\!)$ spanned by the homogeneous formal Laurent series. The cohomology $\rHP_*(\rC):=\rH^*(\rCP(\rC)\!)$ will be called the {\it periodic homology} of $\rC$. Obviously,
$\rHP_*(\rC)\stackrel{u\cdot}\simeq \rHP_{*+2}(\rC)$ (whence ``periodic'') which allows us to view the subscript as taking only two values - $\ev/\od$ (``even/odd'').

Let us denote by $(\Lm,d_\Lambda)$ the unital dg algebra obtained by adjoining to $\Lmo$ generators 
\begin{eqnarray*}
\ee_1,\ldots,\,\ee_{m}\quad  \deg \ee_i=0, \qquad \EE_1,\ldots,\,\EE_{m}\quad  \deg \EE_i=-2 \quad\forall i
\end{eqnarray*}
with the following properties\footnote{$u$ in these formulas should be viewed as a bookkeeping device.}:
\begin{eqnarray}\label{-1}
d_\Lambda(\ee_i+u\EE_i)=[\ee_i+u\EE_i, u\BB_0]-u\BB_i,\quad [\ee_i+u\EE_i, \BB_j]=0 \quad \forall\, i,j=1,\ldots, m.
\end{eqnarray}

Given a (left) dg $\Lm$-module $\rC=(\rrC^*, d_{\rC})$,  $\rHP_*(\rC)$ will stand for the periodic homology of the underlying $\Lmo$-module. The action of the generators $\ee_i$ and $\EE_i$ allows one to extend the $\C\yyy$-module structure on $\rHP_*(\rC)$ to a $\DNC$-module one where $\DNC \supset\C\yyy$ is the algebra obtained by adjoining to $\C\yyy$ the generators $\widetilde{\partial}_{i}$ ($i=1,\ldots, m$) satisfying the relations
\[
\widetilde{\partial}_{i}\cdot \tau_j=\tau_j\cdot\widetilde{\partial}_{i}+\delta_{ij}.
\]
Namely, it follows from (\ref{-1}) that the (degree 0) operators 
$
\widetilde{\partial}_{i}:=\frac{\partial}{\partial \tau_i}-\ee_i-u\EE_i
$
on $\rrC^*\yyy\uu$ commute with the differential 
$d_\rC+u\BB_0+u\sum_i\tau_i\BB_i$ and therefore descent to the periodic homology.  
Obviously, the above-mentioned periodicity  $\rHP_*(\rC)\simeq \rHP_{*+2}(\rC)$ respects the $\DNC$-action.   

To summarize, we have a functor $\rHP_*$ from the category of dg $\Lm$-modules to the category of $\ZZ$-graded $\DNC$-modules. This functor, in fact, extends to a larger category. Namely, let $\rC$ and $\rC'$ be two dg $\Lm$-modules and $\phi: \rC\to\rC'$ be a morphism of the underlying dg $\Lmo$-modules, i.~e.
\begin{eqnarray}\label{dif}
(\ee_i+u\EE_i)\phi-\phi(\ee_i+u\EE_i)
\end{eqnarray}
is not necessarily equal to 0. It follows from (\ref{-1}) that (\ref{dif}), when extended by linearity to a map $\rrC^*\yyy\uu\to\rrC'^*\yyy\uu$, is a morphism of complexes $\rCP(\rC)\to\rCP(\rC')$. Let us say that $\phi$ is a {\it weak} morphism of dg $\Lm$-modules if 
\[
(\ee_i+u\EE_i)\phi-\phi(\ee_i+u\EE_i): \rCP(\rC)\to\rCP(\rC')
\]
is null-homotopic. Obviously, in this case $\rHP_*(\phi)$ is still a morphism of $\DNC$-modules.

Let us conclude by recalling one useful property of the functor $\rHP_*$ (cf. \cite[Prop. 2.4]{GJ}).
\begin{proposition}\label{GJ}
If $\phi: \rC\to\rC'$ is a weak quasi-isomorphism of dg $\Lm$-modules then $\rHP_*(\phi)$ is an isomorphism of $\DNC$-modules. 
\end{proposition}

\subsection{}\label{HC}  
Let us introduce canonical dg $\Lm$-modules associated with dg $\A$-categories.  In what follows, $\cA$, $\cA_\circ$, and $\ttt$ have the same meaning as in Section \ref{dgac}.

The underlying dg $\Lo$-module ($=$ mixed complex)  $(\rrC^*, d_{\rC}, \BB_0)$ of the canonical dg $\Lm$-module associated with $(\cA, \ttt)$ is the standard cyclic mixed complex $(\rrCn_*(\cA), \bb, \B)$ of the underlying dg category (in particular,  $(\rrCn_*(\cA), \bb)$ is the normalized Hochschild complex). Let us start by recalling its description.

Consider the graded space 
\begin{eqnarray*}
\bigoplus_{n\geq0}\bigoplus_{\{\obj_0,\ldots, \,\obj_n\}\subset\cA} 
{\cA}^*(\obj_n,\obj_0)\otimes {\cA}^*(\obj_{n-1},\obj_n)[1]\otimes\ldots\otimes
{\cA}^*(\obj_0,\obj_1)[1].
\end{eqnarray*}
Let us write the elements of the tensor product in the right-hand side as $a_0[a_n|\ldots |a_1]$ ($a_0\in{\cA}(\obj_n,\obj_0)$, $a_p\in{\cA}(\obj_{p-1},\obj_p)$). Then $\rrCn_*(\cA)$ is defined as the quotient of the above space by the subspace spanned by all the tensors of the form $\ldots[\ldots|{\rm id}_{\obj}|\ldots]$ ($\obj\in\cA$). The differential $\bb$
is defined by the formula 
\begin{eqnarray*}
\bb(a_0[a_n|\ldots |a_1])=d_{\cA}(a_0)[a_n|\ldots |a_1]-\sum\limits_{p=1}^n(-1)^{\eta_{p+1}}a_0[a_n|\ldots |d_{\cA}(a_p)|\ldots|a_1]+\\
+(-1)^{|a_0|}a_0a_n[a_{n-1}|\ldots |a_1]+\sum\limits_{p=2}^{n}(-1)^{\eta_p}a_0[a_n|\ldots
|a_pa_{p-1}|\ldots|a_1]
-(-1)^{\eta_2(|a_1|+1)}a_1a_0[a_n|\ldots |a_{2}].
\end{eqnarray*}
where $\eta_p:=|a_0|+|a_p|+\ldots+|a_{n}|+n-p+1$ ($p=1,\ldots, n$) and $\eta_{n+1}:=|a_0|$.
The second differential $\B$ is given by the formula
\begin{eqnarray*}
\B(a_0[a_n|\ldots |a_1])=\sum_{p=0}^{n}(-1)^{\eta_1(\eta_{p+1}+1)}{\rm id}_{\obj_p}[a_p|\ldots|a_1|a_0|\ldots|a_{p+1}].
\end{eqnarray*}

Let us describe now the rest of the canonical dg $\Lm$-module structure on $\rrCn_*(\cA)$. It is this part --  the action of the generators $\BB_i$, $\ee_i$, and $\EE_i$ ($i={1,\ldots,m}$)-- that involves the collection $\ttt=(t_1,\ldots,t_m)\subset\rZ(\cA)$. We set
\[
\BB_i=\BB(t_i),\quad \ee_i=\ee(t_i),\quad \EE_i=\EE(t_i)
\]
where
\begin{eqnarray*}
\BB(t)(a_0[a_n|\ldots |a_1])=-\sum\limits_{p=1}^{n+1}(-1)^{\eta_p}a_0[a_n|\ldots|a_p|t_{\obj_{p-1}}|a_{p-1}|\ldots |a_1],
\end{eqnarray*}
\begin{eqnarray*}
\ee(t)(a_0[a_n|\ldots |a_1])=t_{\obj_0}a_0[a_n|\ldots|a_1],
\end{eqnarray*}
and 
\begin{eqnarray*}
\EE(t)(a_0[a_n|\ldots|a_1])=
\sum_{p=1}^{n+1}\sum_{l=0}^{p-1}(-1)^{\eta_p+(\eta_1-1)\eta_{l+1}}\id_{\obj_l}[a_l|a_{l-1}|\ldots|a_0|a_p|t_{\obj_{p-1}}|a_{p-1}|\ldots |a_{l+1}].
\end{eqnarray*}

\begin{proposition} These operators do define a dg $\Lm$-module structure on $(\rrCn_*(\cA), b)$.
\end{proposition}
\noindent{\bf Proof.} Cf. the proof of \cite[Prop. 3.1, 3.3]{Shk2}. \hfill$\blacksquare$

We will denote the above dg $\Lm$-module and its periodic homology by $\rC(\cA,\ttt)$  and $\rHP_*(\cA,\ttt)$.

The definitions generalize easily to the case of dg $\A$-pairs. Namely, observe that a dg $\A$-functor $F:(\cA,\ttt)\to (\cA',\ttt')$ induces morphisms $\rC(F): \rC(\cA,\ttt)\to\rC(\cA',\ttt').$
In particular, given a full dg subcategory $\cA_\circ$ of $\cA$ one has a morphism $\rC(\cA_\circ,\ttt_\circ)\to\rC(\cA,\ttt)$ where $\ttt_\circ=(t_{1\,\circ},\ldots,t_{m\,\circ})$ stand for the restrictions of the natural transformations onto the subcategory. Then the dg $\Lm$-module associated with $(\cA,\cA_\circ,\ttt)$ is defined by
\begin{equation*}
\rC(\cA,\cA_\circ,\ttt):={\rm Cone}\left(\rC(\cA_\circ,\ttt_\circ)\to\rC(\cA,\ttt)\right).
\end{equation*} 

\subsection{}

Unlike  generic dg $\Lm$-modules, the dg $\Lm$-modules associated with dg $\A$-categories have the following property:
\begin{proposition}\label{wa} The $\DNC$-actions on $\rHP_*(\cA,\ttt)$ and $\rHP_*(\cA,\cA_\circ,\ttt)$ factor through the canonical homomorphism $\DNC\to \Dm$ where $\Dm$ denotes the $m$-th Weyl algebra.
\end{proposition}
\noindent{\bf Proof.} As shown in \cite{Get} (cf. also the proof of  \cite[Prop. 3.5]{Shk2}), there is a canonical degree $-1$ operator $H(t_i,t_j)$ on $\rrCn_*(\cA)\uu$ such that 
$
[\widetilde{\partial}_{i}, \widetilde{\partial}_{j}]=[\bb+u\B+u\sum \tau_i\BB(t_i), H(t_i,t_j)]$.\hfill$\blacksquare$

\medskip
We will also need
\begin{proposition}\label{inv}$\,$

\begin{enumerate}
\item[i)] If a dg $\A$-functor $F:(\cA,\ttt)\to (\cA',\ttt')$ is a dg Morita equivalence of the underlying dg categories then $\rC(F)$ is a quasi-isomorphism (consequently, $\rHP_*(F)$ is an  isomorphism of $\Dm$-modules).
\item[ii)] If a dg $\A$-functor $F:(\cA,\cA_\circ,\ttt)\to (\cA',\cA'_\circ,\ttt')$ of {\rm  localization $\A$-pairs} induces equivalence up to factors of the triangulated categories 
$\rH^0(\cA)/\rH^0(\cA_\circ)\simeq\rH^0(\cA')/\rH^0(\cA'_\circ)
$
then $\rC(F)$ is a quasi-isomorphism (consequently, $\rHP_*(F)$ is an  isomorphism of $\Dm$-modules).
\end{enumerate}
\end{proposition}
\noindent{\bf Proof.} In both cases, the claim is that $F$ induces a quasi-isomorphism on the level of the Hochschild complexes which is already known; cf. \cite[Thm. 5.2]{K2}, \cite[Thm. 2.4]{K1}. \hfill$\blacksquare$

\section{Twisted cyclic homology in the geometric setting}\label{mt} 

We use the notation and keep the assumptions of Section \ref{ex}. In particular, $\cX$ will stand for a smooth variety and $\ww=(f_1,\ldots,f_m)$ for a collection of regular functions on $\cX$.

\subsection{}
As an immediate corollary of Proposition \ref{inv} we obtain
\begin{theorem}\label{mt0} If there exists a complex $\cK$ such that (\ref{phi}) is an equivalence of triangulated categories then the $\Dm$-modules $\rHP_*(\cP(\cX,\ww)\!)$ and $\rHP_*(\cP(\cX',\ww')\!)$ are isomorphic.
\end{theorem}
\noindent{\bf Proof.} As we know (cf. (\ref{Phi})\!), the equivalence lifts to a dg $\A$-functor $\cP^{\rm \,str}(\cX,\ww)\to\cP^{\rm \,flas}(\cX',\ww')$. Since the resulting dg $\A$-functor  satisfies the assumption of  Proposition \ref{inv} ii),  one has an isomorphism of $\Dm$-modules $\rHP_*(\cP^{\rm \,str}(\cX,\ww)\!)\simeq\rHP_*(\cP^{\rm \,flas}(\cX',\ww')\!)$.  
Furthermore,  the embeddings 
\[
\cP^{\rm \,flas}(\cX,\ww)\subset \cP(\cX,\ww)\supset\cP^{\rm \,str}(\cX,\ww) 
\]
also satisfy the condition of Proposition \ref{inv} ii). Thus, 
the $\Dm$-modules  
\[
\rHP_*(\cP(\cX,\ww)\!),\quad
\rHP_*(\cP^{\rm \,str}(\cX,\ww)\!), \quad \rHP_*(\cP^{\rm \,flas}(\cX,\ww)\!),  
\]
are isomorphic.  
 \hfill$\blacksquare$

It may be useful to keep in mind that, in fact, all the $\Dm$-modules  
\[
\rHP_*(\cP(\cX,\ww)\!),\,
\rHP_*(\cP^{\rm \,str}(\cX,\ww)\!),\, \rHP_*(\cP^{\rm \,inj}(\cX,\ww)\!), \,\rHP_*(\cP^{\rm \,flas}(\cX,\ww)\!), \, \rHP_*({\rm End}(\cE), \ww)
\]
are isomorphic since the embeddings 
\[
({\rm End}(\cE), \ww)\subset \cP^{\rm \,inj}(\cX,\ww)\quad \text{and}\quad  \cP^{\rm \,inj}(\cX,\ww)\subset \cP^{\rm \,flas}(\cX,\ww)\subset \cP(\cX,\ww)\supset\cP^{\rm \,str}(\cX,\ww) 
\]
satisfy the conditions of Proposition \ref{inv} i) and ii), respectively.

\subsection{}

Consider the following multi-variable version of the complex we discussed in the Introduction:
\[
(\Omega^*_\cX\yyy,\, d-\sum_{i=1}^m\tau_i df_i).
\] 
It is a complex of  $\Dm$-modules where the action of $\Dm$ is given by the operators of multiplication with the $\tau_i$'s and the operators $\frac{\partial}{\partial \tau_i}-f_i$ (obviously, all these operators commute with the differential). The $\Dm$-action induces one on the hypercohomology $\mathbb{H}^{k}(\Omega^*_\cX\yyy,\, d-\sum \tau_i df_i)$ which can be seen as follows. The hypercohomology can be identified with the cohomology of the complex
\begin{eqnarray}\label{9}
\left(\bigoplus_{p+q=*}\Gamma(G^p(\Omega^q_\cX)\!)\yyy, \,\,\Gamma(G(d)\!)+(-1)^{\star}\Gamma({d_G})-\sum \tau_i \Gamma(G(df_i)\!)\right)
\end{eqnarray}
where $\cF\to G(\cF)=(G^*(\cF), {d_G})$ denotes the (functorial) Godement resolution, $\Gamma$ the global section functor, and $\star|_{\Gamma(G^p(\Omega_\cX^q)\!)}:=q\cdot {\rm id}$. The $\Dm$-action on (\ref{9}) is given by the operators of multiplication with the $\tau_i$'s and the operators $\frac{\partial}{\partial \tau_i}-\Gamma(G(f_i)\!)$.

\begin{theorem}\label{mt1} There is an isomorphism of $\Dm$-modules
\begin{eqnarray}\label{mr1}
\rHP_{\ev/\od}(\cP(\cX,\ww)\!)\simeq \bigoplus_{k\,\text{\rm even/odd}} \mathbb{H}^{k}(\Omega^*_\cX\yyy,\, d-\sum_{i=1}^m\tau_i df_i).
\end{eqnarray}
\end{theorem}
\noindent{\bf Proof.} The proof of this theorem relies on (and is an extension of) the proof of Theorem 5.2 in \cite{K0} which in combination with results in \cite{We} establishes the claim in the case $m=0$. Therefore we will present a sketch of the proof -- modulo results and some technical details that can be found in \cite{K0}. 
A complementary reference is \cite{Ef} where a similar argument is employed for calculating the cyclic mixed complexes of the categories of matrix factorizations associated with $(\cX,f)$.  

\subsubsection{}\label{ss1} The proof of the theorem is based on the reduction to the case when $\cX$ is affine. This requires working with (pre)sheaves of dg $\Lm$-modules on $\cX$.
The case of sheaves of ordinary mixed complexes ($m=0$) is discussed in detail in \cite{We,K0,Ef}, and all the basic definitions and constructions generalize to an arbitrary $m$ in a straightforward manner. In particular, presheaves of dg $\Lm$-modules can be sheafified and there is a version of the derived global section functor $\rg$ in this setting. Let us discuss the latter notion.

Recall that $\cF\to G(\cF)=(G^*(\cF), d_G)$ stands for the Godement resolution\footnote{A comprehensive discussion of the Godement resolution and its application to computing the total right derived global section functor can be found in \cite{RGR} (cf. especially Section 5.2).} of sheaves on $\cX$.  Thanks to its functoriality one has the following endofunctor on the category of sheaves of dg $\Lm$-modules:
\[
\rCu=(\rrCu^*, d_\rrCu) \to G(\rCu):=(\prod_{p+q=*}G^p(\rrCu^{q}), G(d_\rrCu)+(-1)^\star d_G).
\]
Then one defines $\rg(\rCu):=\Gamma(G(\rCu)\!)$ where the global section functor $\Gamma$ in the right-hand side is applied component-wise. Thus, $\rg$ takes values in the category of ordinary dg $\Lm$-modules.

The most important property of $\rg$ is that  it reflects quasi-isomorphisms. By a quasi-isomorphism between two sheaves $\rCu$ and $\rCu'$ of dg $\Lm$-modules we understand a morphism $\rCu\to\rCu'$ that induces a quasi-isomorphism of the underlying complexes. (In turn, the latter means that we have quasi-isomorphisms stalk-wise: $\rrCu^*_x\simeq\rrCu^{\prime *}_x$ for all $x\in\cX$. In practice, it is enough to check that for every open affine subset $U\subset\cX$ the morphisms $\rrCu^*(U)\to\rrCu^{\prime *}(U)$ are quasi-isomorphisms since the open affine subsets form a base of the Zariski topology.)

\subsubsection{} A key observation is that $\rM(\cP(\cX,\ww)\!)$ is quasi-isomorphic to the derived global sections of a sheaf of dg $\Lm$-modules. 
Namely, the assignment
\begin{eqnarray*}
U\mapsto \cM(\cP(\cX,\ww)\!)(U):=\rM(\cP\, U,\cP_{\circ} U,\ww|_U) \quad (U\subset\cX\,\, \text{is open})
\end{eqnarray*}
 is a presheaf of dg $\Lm$-modules, and  the sheaf in question is its sheafification, $\cM^{\rm sh}(\cP(\cX,\ww)\!)$. More precisely, the claim is that the composition of the natural morphisms
\begin{eqnarray*}
\rM(\cP(\cX,\ww)\!)=\Gamma(\cM(\cP(\cX,\ww)\!)\!)\to\Gamma(\cM^{\rm sh}(\cP(\cX,\ww)\!)\!)\to\rg(\cM^{\rm sh}(\cP(\cX,\ww)\!)\!)
\end{eqnarray*}
is a quasi-isomorphism. Note that the claim is about the underlying complexes and its validity has nothing to do with $\ww$. For the underlying complexes the result is already known to hold; cf.   Lemma 4.1,  Diagram 4.2, and Theorem 5.2 in \cite{K0}.

As a consequence, 
\begin{eqnarray}\label{r2}
\rHP_*(\rM(\cP(\cX,\ww)\!)\!)\simeq\rHP_*(\rg(\cM^{\rm sh}(\cP(\cX,\ww)\!)\!)\!).
\end{eqnarray}  In particular, by Proposition \ref{wa} the $\DNC$-action on $\rHP_*(\rg(\cM^{\rm sh}(\cP(\cX,\ww)\!)\!)\!)$  factors through the homomorphism $\DNC\to \Dm$.

\subsubsection{}\label{ss2} Let us denote by $\Omega'(\cX,\ww)$ the sheaf of dg $\Lm$-modules on $\cX$ whose underlying complex of sheaves is $(\rrCu^*, d_\rrCu):=(\Omega^{-*}_\cX, 0)$
(note the reversal of the grading) and the action of $\Lm$ is given by
\begin{eqnarray*} 
\B:=d, \quad \BB_i:=-df_i,\quad \ee_i:=f_i, \quad \EE_i:=\frac12 df_i\cdot d \quad (i=1,\ldots, m)
\end{eqnarray*}
We claim that $\cM^{\rm sh}(\cP(\cX,\ww)\!)$ and $\Omega'(\cX,\ww)$ are quasi-isomorphic. 

To prove this, consider the presheaf of dg $\Lm$-modules $\cM(\cO_\cX,\ww)$ defined as follows. Associated with an open subset $U\subset\cX$ is the dg $\A$-algebra $(\cO_\cX(U),\ww|_U)$ (with the trivial grading and differential), and we set
$
 \cM(\cO_\cX,\ww)(U) :=\rM(\cO_\cX(U),\ww|_U). 
$
According to \cite[Sect. 3.4, 4.1]{K0}, there is a quasi-isomorphism from the underlying complex of  $\cM^{\rm sh}(\cO_\cX,\ww)$ to the underlying complex of $\cM^{\rm sh}(\cP(\cX,\ww)\!)$.  This quasi-isomorphism is easily checked to preserve the $\Lm$-action (essentially, this is due to the fact that the morphism comes from dg $\A$-functors). 

Let us construct a quasi-isomorphism $\cM^{\rm sh}(\cO_\cX,\ww)\to\Omega'(\cX,\ww)$. Let $U\subset\cX$ be any open subset and consider
the classical Hochschild-Kostant-Rosenberg map $\varepsilon_U: \rrCn_*(\cO_\cX(U)\!)\to \Omega^{-*}_\cX(U)$  given by 
\[
\varepsilon(\phi_0[\phi_1|\ldots |\phi_n])=\frac{1}{n!}\phi_0d\phi_1\wedge\ldots \wedge d\phi_n, \quad \phi_i\in\cO_\cX(U).
\]
It is well known that $\varepsilon$ gives rise to a morphism of the corresponding dg $\Lo$-modules and that this morphism is a quasi-isomorphism when $U$ is affine. That $\varepsilon$ preserves the rest of the $\Lm$-action was checked in \cite[Lem. 4.1]{Shk2} and \cite[Sect. D.5]{Shk1}. 

As a consequence, 
\begin{eqnarray}\label{r3}
\rHP_*(\rg(\cM^{\rm sh}(\cP(\cX,\ww)\!)\!)\!)\simeq \rHP_*(\rg(\Omega'(\cX,\ww)\!)\!)
\end{eqnarray}  In particular, the $\DNC$-action on the right-hand side factors through the homomorphism $\DNC\to \Dm$. 

\subsubsection{} Let $\Omega(\cX,\ww)$ be the sheaf of dg $\Lm$-modules on $\cX$ whose underlying complex of sheaves, the $\BB$- and the $\ee$-operators are exactly the same as for $\Omega'(\cX,\ww)$ but $\EE_i=0$ for all $i$. We claim that there is an isomorphism of $\DNC$-modules (which are, in fact, $\Dm$-modules)
\begin{eqnarray}\label{r4}
\rHP_{*}(\rg(\Omega'(\cX,\ww)\!)\!)\simeq\rHP_{*}(\rg(\Omega(\cX,\ww)\!)\!).
\end{eqnarray}
More precisely, we claim that the identity endomorphism of the underlying complex defines a  weak quasi-isomorphism of dg $\Lm$-modules $\rg(\Omega'(\cX,\ww)\!)\to\rg(\Omega(\cX,\ww)\!)$. Indeed, let \[(K^*, d_K):=\rCP(\rg(\Omega'(\cX,\ww)\!)\!)=\rCP(\rg(\Omega(\cX,\ww)\!)\!),\] i.~e.
\begin{eqnarray*}
K^*=\bigoplus\limits_{p-q+2l=*} \Gamma(G^p(\Omega_\cX^q)\!)\yyy u^l, \quad d_K=(-1)^\star\Gamma({d_G})+u\Gamma(G(d)\!)-u\sum \tau_i \Gamma(G(df_i)\!).
\end{eqnarray*}
The operators $\ee_i+u\EE_i$ in the two cases are equal to
\begin{equation*}
\Gamma(G(f_i)\!)+\frac{u}2 \Gamma(G(df_i\cdot d)\!)\quad \text{and}\quad \Gamma(G(f_i)\!),
\end{equation*}
respectively. We need to show that the difference $\frac{u}2 \Gamma(G(df_i\cdot d)\!)$,  viewed as an endomorphism of $(K^*, d_K)$, is homotopic to 0. It is easy to check that
\[
u \Gamma(G(df_i\cdot d)\!)=\Gamma(G(f_i\cdot d)\!)\cdot d_K+d_K\cdot\Gamma(G(f_i\cdot d)\!).
\]

\subsubsection{} Let $(L^*, d_L)$ denote the complex (\ref{9}) computing $\mathbb{H}^*(\Omega^*_\cX\yyy,\, d-\sum \tau_i df_i)$ and $(K^*, d_K)$ be as above.
It is easy to see that the map
\[
u^{-\star}: \bigoplus\limits_{p,q}\Gamma(G^p(\Omega_\cX^q)\!)\to\bigoplus\limits_{p,q}\Gamma(G^p(\Omega_\cX^q)\!)\otimes\C[\!u\!]
\]
(recall that $\star|_{\Gamma(G^p(\Omega_\cX^q)\!)}=q\cdot {\rm id}$) extends to an isomorphism of complexes $(K^*, d_K)\to (L^*\uu, d_L)$ which, in addition, commutes with the operators $\frac{\partial}{\partial \tau_i}-\Gamma(G(f_i)\!)$. This implies
\[
\rHP_{\ev/\od}(\rg(\Omega(\cX,\ww)\!)\!)\simeq \bigoplus_{k\,\text{\rm even/odd}} \mathbb{H}^{k}(\Omega^*_\cX\yyy,\, d-\sum \tau_i df_i).
\]
which together with (\ref{r2}), (\ref{r3}) and (\ref{r4}) completes the proof of Theorem \ref{mt1}. \hfill$\blacksquare$

\medskip

Theorems \ref{mt0} and \ref{mt1} yield

\begin{corollary}\label{cor} If there exists a complex $\cK$ such that (\ref{phi}) is an equivalence of triangulated categories then there is an isomorphism of $\Dm$-modules
\begin{eqnarray*}
\bigoplus_{k\,\text{\rm even/odd}} \mathbb{H}^{k}(\Omega^*_\cX\yyy,\, d-\sum_{i=1}^m\tau_i df_i)\simeq \bigoplus_{k\,\text{\rm even/odd}} \mathbb{H}^{k}(\Omega^*_{\cX'}\yyy,\, d-\sum_{i=1}^m\tau_i df'_i).
\end{eqnarray*}
\end{corollary}

\end{document}